\documentclass[11pt,twoside]{article}
\usepackage{amssymb,amsmath,mathrsfs,txfonts,graphicx,color}
\graphicspath{{figs/}}

\begin{document}

\def\pd#1#2{\frac{\partial#1}{\partial#2}}
\def\dfrac{\displaystyle\frac}
\let\oldsection\section
\renewcommand\section{\setcounter{equation}{0}\oldsection}
\renewcommand\thesection{\arabic{section}}
\renewcommand\theequation{\thesection.\arabic{equation}}

\newtheorem{theorem}{\indent Theorem}[section]
\newtheorem{lemma}{\indent Lemma}[section]
\newtheorem{proposition}{\indent Proposition}[section]
\newtheorem{definition}{\indent Definition}[section]
\newtheorem{remark}{\indent Remark}[section]
\newtheorem{corollary}{\indent Corollary}[section]

\title{Theoretical and numerical studies on global stability of traveling waves with oscillations
for time-delayed nonlocal dispersion equations
}
\author{Tianyuan Xu$^a$, Shanming Ji$^{b,}$\thanks{Corresponding author.
{\it E-mail}:jism@scut.edu.cn},
\ Rui Huang$^a$, Ming Mei$^{c,d}$, Jingxue Yin$^a$
\\
\\
{ \small \it $^a$School of Mathematical Sciences, South China Normal University}
\\
{ \small \it Guangzhou, Guangdong, 510631, P.~R.~China}
\\
{ \small \it $^b$School of Mathematics, South China University of Technology}
\\
{ \small \it Guangzhou, Guangdong, 510641, P.~R.~China}
\\
{ \small \it $^c$Department of Mathematics, Champlain College Saint-Lambert}
\\
{ \small \it Quebec,  J4P 3P2, Canada, and}
\\
{ \small \it $^d$Department of Mathematics and Statistics, McGill University}
\\
{ \small \it Montreal, Quebec,   H3A 2K6, Canada}
\\[3pt]
}
\date{}

\maketitle

\begin{abstract}
This paper is concerned with the global stability of non-critical/critical traveling waves with oscillations
for time-delayed nonlocal dispersion equations.
We first theoretically prove that all traveling waves,
especially the critical oscillatory traveling waves,
are globally stable in a certain weighted space, where the convergence rates to the non-critical oscillatory traveling waves are
time-exponential, and the convergence to the critical oscillatory traveling waves are
time-algebraic. Both of the rates are optimal.
The approach adopted is the weighted energy method
with the fundamental solution theory for time-delayed equations. Secondly, we carry out
numerical computations in different cases, which also confirm our theoretical results.
Because of oscillations of the solutions and nonlocality of the equation, the numerical results obtained by the
regular finite difference scheme are not stable, even worse to be blow-up. In order to overcome these obstacles,
we propose a new  finite difference scheme by adding artificial viscosities to both sides of the equation, and
obtain the desired  numerical  results.

\end{abstract}



\section{Introduction}

In this paper, we  consider the global stability of critical
oscillatory traveling waves for a class
of nonlocal dispersion equations with time-delay
\begin{equation} \label{eq-main}
\begin{cases}
\displaystyle
\pd{v}{t}-D(J*v-v)+d(v)=K*b(v(t-r,\cdot)), \quad &x\in\mathbb R, ~t>0,\\
v(s,x)=v_0(s,x), \quad &x\in\mathbb R, ~s\in[-r,0],
\end{cases}
\end{equation}
where the initial value satisfies
$$\lim_{x\to\pm\infty}v_0(s,x)=v_\pm, ~\text{uniformly in}~ s\in[-r,0].$$
This model represents the spatial dynamics of a single-species population
with age-structure and nonlocal diffusion such as the Australian blowflies population
distribution \cite{Gourley1,Gourley2}.
Here $v(t, x)$ denotes the total mature population of the species,
the function $d(v)$ and $b(v)$ are the death and birth rates of the mature
population respectively,
$J(x)$ and $K(x)$ are non-negative, unit and
symmetric kernels,
$J(x)$ is the probability distribution of rates of dispersal over distance $x$.
Then $J*v(x)$ is the rate at which individuals are arriving at position $x$ from all other locations,
and $v(x)=\int_\mathbb RJ(x,y)v(x)dy$ is the rate at which they are leaving location $x$ to travel to all
other sites. Therefore, the expression $D(J*v-v)$ is the nonlocal dispersion due to long range dispersion mechanisms \cite{Fife1,Hutson},
where the coefficient $D>0$ is the spatial diffusion rate.

The advantages of the nonlocal process governed by integral process over the
classical dispersal process modelled by Laplacian lie in
the fact that the nonlocal one accounts for interaction between individual in both short and long ranges, while the classical one accounts for only local interactions between the neighbor individuals.
Moreover,
the nonlocal operator for the initial value problem is not a smoothing operator.
Discontinuities in the initial data are retained \cite{Fife1}.
And
the spatial decay rates of the traveling waves at  infinity are different in the
local and nonlocal cases \cite{Wang02}.



From the classical Nicholson's blow flies model \cite{Gurney-Blythe} with the birth rate function $b(v)=pve^{-av}$ for $p>0$ and $a>0$ and the death rate
function $d(v)=\delta v$ for $\delta>0$, and the Mackey-Glass model \cite{Mackey-Glass} with $b(v)=\frac{v}{1+av^q}$ for $a>0$ and $q>1$ and
$d(v)=dv$ for $d>0$,   throughout this paper, we assume the birth rate function, the death rate function, and the kernels to be:

(H1) There exist two constant equilibria of \eqref{eq-main}:
$v_-=0$ is unstable and $v_+>0$ is stable, namely,
$d(0)=b(0)=0$, $d(v_+)=b(v_+)$, $0\le d'(0)<b'(0)$ and $d'(v_+)>b'(v_+)$;

(H2) Both $d(s)$ and $b(s)$ are non-negative, $C^2$-smooth functions with $d'(s)\ge d'(0)\ge0$, $|b'(s)|\le b'(0)$ for $s\in[0,+\infty)$,
but $b(s)$ is non-monotone;

(H3) Both kernels $J(x)$ and $K(x)$ are
nonnegative, symmetric and unit,
\begin{align*}
J(x)\ge0,\quad J(-x)=J(x), \quad
\int_\mathbb R J(x)dx=1, \\
K(x)\ge0,\quad K(-x)=K(x),\quad
\int_\mathbb R K(x)dx=1,
\end{align*}
and satisfy
$$
\int_\mathbb R |x|J(x)e^{-\eta x}dx<\infty
\quad
\text{and~}
\int_\mathbb R|x|K(x)e^{-\eta x}dx<\infty
\quad\text{for any}\quad \eta>0.
$$

(H4) The Fourier transform of $J(x)$, denoted by $\hat J(\xi)$, satisfies that
$\hat J(\xi)=1-\kappa|\xi|^\alpha+o(|\xi|^\alpha)$ as $\xi\to0$
with $\alpha\in(0,2]$ and $\kappa>0$,
and $1-\hat J(\xi)\ge \omega(r)$ for all $|\xi|\ge r$ and any $r>0$
with some positive function $\omega(r)>0$.

A traveling wavefront of \eqref{eq-main} is a special solution of the form $u(t, x)=\phi(x+ct)$, where $c$ is the wave speed.
The existence and uniqueness (up to a shift) of traveling waves for the
equation \eqref{eq-main} were proved in  \cite{HuangMeiWang,ZhangNATMA11,ZhangAA}.
The main purpose of this paper is
to study the global stability of traveling wavefronts $\phi(x+ct)$ of \eqref{eq-main},
especially the case of the critical wave $\phi(x+c^*t)$.
Here the number $c^*$ is called the
critical speed (or the minimum speed) in the sense that a traveling wave
exists if $c\ge c^*$, while no traveling wave $\phi(x+c^*t)$ exists if $0<c< c^*$.
Let $\phi(x+ct)=\phi(\xi)$ be any given monotone or non-monotone
traveling waves for \eqref{eq-main} with wave speed $c\ge c^*$
connecting the two steady equilibria $v_{\pm}$, namely,
\begin{equation} \label{eq-tw}
\begin{cases}
\displaystyle
c\phi'(\xi)-D\Big(\int_{\mathbb R}J(y)\phi(\xi-y)dy-\phi(\xi)\Big)
+d(\phi(\xi)) \\
\displaystyle
\qquad\qquad=\int_{\mathbb R}K(y)b(\phi(\xi-y-cr))dy, \qquad \xi\in\mathbb R,\\
\phi(\pm\infty)=v_{\pm}, \quad \phi(\xi)\ge0, \qquad \xi\in\mathbb R,
\end{cases}
\end{equation}
where $\xi=x+ct$, ${}'=\pd{}{\xi}$.
As summarized in \cite{HuangMeiWang},
we obtain the
following characteristic equation for the pair of $(c, \lambda)$:
\begin{equation}
c\lambda-D\int_\mathbb R J(y)e^{-\lambda y}dy+D+d'(0)=
b'(0)e^{-\lambda cr}\int_\mathbb R K(y)e^{-\lambda y}dy.
\end{equation}
The critical speed $c^*$ is uniquely determined by
\begin{align}\label{eq-characteristic1}
b'(0)e^{-\lambda_* cr}\int_\mathbb R K(y)e^{-\lambda_* y}dy
=c^*\lambda_*-\int_\mathbb R J(y)e^{-\lambda_*y}dy+D+d'(0),\\\label{eq-characteristic2}
-b'(0)e^{-\lambda_*c^*r}\int_\mathbb R
(y+c^*r)K(y)e^{-\lambda_*y}dy
=c^*+D\int_\mathbb R yJ(y)e^{-\lambda_*y}dy.
\end{align}
and when $c>c^*$, there exist two numbers $\lambda_2>\lambda_1>0$ such that
\begin{equation}\label{eq-lambda12}
b'(0)e^{-\lambda_i cr}\int_\mathbb R K(y)e^{-\lambda_i y}dy
=c^*\lambda_i-\int_\mathbb R J(y)e^{-\lambda_iy}dy+D+d'(0)\quad
\text{for}\quad i=1,2.
\end{equation}
and
\begin{equation}
c^*\lambda-\int_\mathbb R J(y)e^{-\lambda y}dy+D+d'(0)>
b'(0)e^{-\lambda cr}\int_\mathbb R K(y)e^{-\lambda y}dy\quad
\text{for}\quad \lambda\in(\lambda_1,\lambda_2).
\end{equation}
As discussed in \cite{Chern}, when $b'(v_+)<0$ and $|b'(v_+)|\le d'(v_+)$, the traveling waves may occur oscillations around $v_+$ for the time-delay
$r>\underline{r}$, where $\underline{r}$, given by
\[
|b'(v+)|re^{d'(v_+)r+1} = 1,
\]
is the critical point for the solution to the delayed ODE
\[
v'(t) + d'(v_+)v(t) = b'(v+)v(t - r)
\]
and based on Hopf-bifurcation analysis, there will be no traveling waves if the time-delay $r\ge \overline{r}$, where  $\overline{r}$ is the
Hopf-bifurcation point:
 \begin{equation}\label{overline-r-1}
\overline r:=\frac{\pi-\arctan(\sqrt{|b'(v_+)|^2-d'(v_+)^2}/d'(v_+))}
{\sqrt{|b'(v_+)|^2-d'(v_+)^2}}.
\end{equation}

{\color{red}}
There have been extensive investigations on the stability of traveling waves for
reaction-diffusion equations with and without time delay \cite{Chern,Fife,Gallay,Ji-Yin-Huang,LinLinSIAM2014,MeiEdinburgh04,
MeiEdinburgh08,MeiWang2011,MeiWong,Schaaf,Shen,Smith,Wang,Wu}.
For the reaction-diffusion equations with time-delay and local dispersal, the first work on the linear stability of the traveling wave for time-delayed
reaction-diffusion equation was given by Schaaf \cite{Schaaf} in 1987 based on spectral analysis.
For the bistable case, Smith and Zhao \cite{Smith} obtained the stability of traveling waves for local equations by the upper-lower solutions method.
Later then, Wang-Li-Ruan \cite{Wang} proved
the existence and globally asymptotic stability of traveling wave fronts for equations with nonlocal delay.

Compared to the rich results for the local reaction-diffusion equations, limited theoretical results exist for the equations with nonlocal dispersion \cite{MeiLinJDE091,MeiLinJDE092,ZhangNATMA11,ZhangMa,ZhangWangJDE12} .
When the birth rate $b(v)$ is monotone,
Pan-Lin-Lin \cite{PanLin} first showed the local stability for the monotone wave when the wave
speed is sufficiently large $c\gg1$ via upper and lower solutions method.
Furthermore,
Huang-Mei-Wang \cite{HuangMeiWang} proved that all noncritical and critical monotone wavefronts are globally stable by  Fourier transform and the weighted energy method.
When the birth rate $b(v)$ is non-monotone,
the equation \eqref{eq-main} losses monotonicity and the solution will be oscillating or even not exists for large time-delay $r$.
In this case, Zhang \cite{ZhangNATMA11} obtained
the existence of traveling waves with $c>c^*$ by introducing two auxiliary nonlocal dispersal equations with quasi-monotonicity.
Zhang-Ma \cite{ZhangMa}
further proved that the traveling waves with sufficiently large speed $c\gg1$ are
locally stable, when the initial perturbation around the wave front is small. 
The asymptotic stability of non-critical oscillatory traveling waves with $c>c^*$ was proved by Huang-Mei-Zhang-Zhnag in \cite{HuangDCDS16}. Note that,  their results still need the assumption of
the small initial perturbations around the waves. The question
whether these oscillatory traveling waves are globally stable for large perturbations is not clear at all.

However, the most interesting cases are for the slower wave speed $c\ge c^*$, especially for
the critical waves.
As mentioned in \cite{Liang,Thieme,ZhangMa}, the critical wave speed coincides with the asymptotic speed of propagation, and it is very important in the biological invasions.
Zhang-Ma \cite{ZhangMa}  established the existence of critical traveling wave solution with $c=c^*$ and proved the number $c^*$ is also the
spreading speed of the corresponding initial value problem with compact support.

The stability for the critical oscillating traveling
waves of the reaction-diffusion equation with nonlocal dispersion and time delay \eqref{eq-main} with $c=c^*$ is open so far as we know.
In fact, the study on stability of critical traveling waves is very limited
and also very challenging.
For local dispersion case,
there are several works on stability of the critical waves of some typical reaction diffusion equations.
In 1979,
Moet \cite{Moet} showed the algebraic stability of the critical waves for the
classical Fisher-KPP equation by the Green function method.
Later on, Gallay \cite{Gallay} improved the  algebraic convergence rate via renormalization
group method.
For the local Nicholson's blowflies equation,
the global stability of critical traveling waves is obtained in \cite{MeiSIAM10},
and the convergence rate to the critical
wave is proved to be algebraic by the Green function method.
Regarding the non-monotone traveling waves, Lin-Lin-Lin-Mei \cite{LinLinSIAM2014} first proved that all non-critical
non-monotone traveling wave are time-exponentially stable when the initial perturbations around the waves are small enough. Furthermore,
Chern-Mei-Yang-Zhang \cite{Chern} proved that all critical non-monotone  traveling waves are locally
stable by the anti-weighted energy method, but, due to the technical reason, there is no convergence rate addressed.
Very recently, Mei-Zhang-Zhang \cite{MeiZhang} developed a new method to prove the global stability of critical oscillating traveling waves.
They based on some key observations for the structure of the govern equations establishing the boundedness estimate for the oscillating solutions.
Inspired by their work, we intend to study the nonlocal dispersion equation with time-delay.

Our main purpose is to study the global stability of the reaction-diffusion equation \eqref{eq-main} with nonlocal dispersion and time delay for all traveling waves, especially the critical traveling waves.
Due to the lack of monotonicity, the bad effect of time delay and the nonlocality,
we have to face some new challenges and  look for a new strategy to solve
the problem.
The main approach adopted is the weighted energy method with some new developments.
We prove that for all oscillatory traveling waves, including the critical traveling waves are globally stable, where the initial perturbations in a certain weighted Sobolev space can be arbitrarily big. The convergence to the non-critical traveling waves with $c>c^*$ is time-exponential, and the convergence to the critical traveling waves with $c=c^*$ is time-algebraic.

We define the uniformly continuous space $C_\text{unif}[-r,T]$ for $0<T\le\infty$,
by
\begin{align*}
C_\text{unif}[-r,T]:=&\{v(t,x)\in C([-r,T]\times \mathbb R)~\text{such that}\\
&\lim_{x\to+\infty}v(t,x)~\text{exists uniformly in}~t\in[-r,T], \text{and}\\
&\lim_{x\to+\infty}v_x(t,x)=0,
\text{uniformly with respect to} \quad t\in[-r,T]\}.
\end{align*}
For $c\ge c^*$, we define the following weight function
\begin{equation} \label{eq-weight}
w(\xi)=
\begin{cases}
e^{-2\lambda\xi}, ~\xi\in\mathbb R,
\qquad &\text{for~} c>c^*, \lambda\in(\lambda_1,\lambda_2),\\
e^{-2\lambda^*\xi}, ~\xi\in\mathbb R,
\qquad &\text{for~} c=c^*,
\end{cases}
\end{equation}
where $\lambda_1$ and $\lambda_2$ are specified in \eqref{eq-lambda12}.
Notice that for $c\ge c^*$,
$\displaystyle\lim_{\xi \to -\infty}w(\xi)=+\infty$
and $\displaystyle\lim_{\xi \to -\infty}w(\xi)=0$,
because $\lambda>0$ and $\lambda_*>0$.
We denote the weighted Sobolev spaces $L_w^1(\mathbb R)$ and $H_w^{1}(\mathbb R)$ by
$$L_w^1(\mathbb R)=\{u; wu\in L^1(\mathbb R)\},$$
and
$$H_w^{1}(\mathbb R)=\{u; \sqrt{w}u, \sqrt{w} u_x\in L^2(\mathbb R)\}.$$

Our main stability theorems are as follows.
\begin{theorem}[Global stability] \label{th-stability}
Assume that (H1)--(H4) hold.
Let $b'(v_+)$ and $r$ satisfy,
either $d'(v_+)\ge|b'(v_+)|$ with arbitrary $r>0$,
or $d'(v_+)<|b'(v_+)|$ with $0<r<\overline r$,
where $\overline r$ is defined in \eqref{overline-r-1}.
Let $\phi(\xi)=\phi(x+ct)$ be any given traveling wave with $c\ge c^*$
and the initial perturbation be
$v_0(s,x)-\phi(x+cs)\in C_\text{unif}[-r,0]\cap C([-r,0];L_w^1(\mathbb R)\cap H_w^1(\mathbb R))$
and $\partial_s(v_0-\phi)\in C([-r,0];L_w^1(\mathbb R)\cap H_w^1(\mathbb R))$.
Then the global solution $v(t,x)$ of \eqref{eq-main} satisfies

i) if $c>c^*$, then
$$\sup_{\mathbb R}|v(t,x)-\phi(x+ct)|\le Ct^{-\frac{1}{\alpha}}e^{-\mu t}$$
for some positive constants $\mu>0$ and $C>0$;

ii) if $c=c^*$, then
$$\sup_{\mathbb R}|v(t,x)-\phi(x+ct)|\le Ct^{-\frac{1}{\alpha}}$$
for some positive constant $C>0$.
\end{theorem}

This paper is organized as follows. In Section 2, we prove our main stability theorem. Then we shall carry out  numerical simulations for Nicholson's blowflies model with nonlocal diffusion in Section 3, which further numerically confirm
our theoretical results.

\section{Proof of the main results}

Now we consider the perturbed solution $v(t,x)$ of \eqref{eq-main}
around any given traveling waves $\phi(x+ct)=\phi(\xi)$ of \eqref{eq-tw}.
Define
\begin{align*}
u(t,\xi):=&v(t,x)-\phi(x+ct)=v(t,\xi-ct)-\phi(\xi), \\
u_0(s,\xi):=&v_0(s,x)-\phi(x+cs)=v_0(s,\xi-cs)-\phi(\xi).
\end{align*}
Then $u(t,\xi)$ satisfies
\begin{equation} \label{eq-pt}
\begin{cases}
\displaystyle
\pd{u}{t}
+c\pd{u}{\xi}-D\Big(\int_{\mathbb R}J(y)u(t,\xi-y)dy-u(t,\xi)\Big)
+P[u](t,\xi) \\
\displaystyle
\qquad\qquad=\int_{\mathbb R}K(y)Q[u](t-r,\xi-y-cr) dy,
\quad \xi\in\mathbb R,~t>0\\
u(s,\xi)=u_0(s,\xi), \quad \xi\in\mathbb R, ~s\in[-r,0],
\end{cases}
\end{equation}
where
\begin{align} \label{eq-PQ}
\begin{cases}
P[u](t,\xi):=d(\phi(\xi)+u(t,\xi))-d(\phi(\xi)),\\
Q[u](t,\xi):=b(\phi(\xi)+u(t,\xi))-b(\phi(\xi)).
\end{cases}
\end{align}

We first show the existence and uniqueness of solution $u(t,\xi)$
to the initial value problem of time-delayed nonlocal dispersion equation
\eqref{eq-pt} in the uniformly continuous space $C_\text{unif}[-r,+\infty)$.

\begin{lemma} \label{le-existence}
Assume (H1)-(H3) hold.
If the initial perturbation $u_0\in C_\text{unif}[-r,0]$,
then the perturbed problem \eqref{eq-pt} admits one unique global solution
$u(t,\xi)$ in $C_\text{unif}[-r,+\infty)$.
\end{lemma}
{\it\bfseries Proof.}
First we solve the problem for $t\in[0,r]$.
Since $t-r\in[-r,0]$ and $u(t-r,\xi)=u_0(t-r,\xi)$, \eqref{eq-pt} is reduced to
\begin{equation} \label{eq-zpt}
\begin{cases}
\displaystyle
\pd{u}{t}
+c\pd{u}{\xi}-D\Big(\int_{\mathbb R}J(y)u(t,\xi-y)dy-u(t,\xi)\Big)
+P[u](t,\xi) \\
\displaystyle
\qquad\qquad=\int_{\mathbb R}K(y)Q[u_0](t-r,\xi-y-cr) dy,
\quad \xi\in\mathbb R,~t>0\\
u(0,\xi)=u_0(0,\xi), \quad \xi\in\mathbb R,
\end{cases}
\end{equation}
Back to the original coordinates, that is,
we make change of variable, $u(t,\xi)=u(t,x+ct)=\bar u(t,x)$, $\xi=x+ct$,
the above problem \eqref{eq-zpt} is equal to
\begin{equation} \label{eq-zptx}
\begin{cases}
\displaystyle
\pd{\bar u}{t}
-D\Big(\int_{\mathbb R}J(y)\bar u(t,x-y)dy-\bar u(t,x)\Big)
+P[\bar u](t,x) \\
\displaystyle
\qquad\qquad=\int_{\mathbb R}K(y)Q[u_0](t-r,x-y) dy,
\quad x\in\mathbb R,~t>0\\
\bar u(0,x)=u_0(0,x), \quad x\in\mathbb R,
\end{cases}
\end{equation}
The existence of solution to \eqref{eq-zptx} follows from the semigroup theory of
the convolution operators. In fact, from the textbook \cite{Rossi}, it is known that
 \eqref{eq-zptx} can be written in the integral form of
\begin{equation}\label{new-1}
u(t,x)=S(t)\ast u_0 - \int^t_0 S(t-\tau)\ast P[u](s) ds + \int^t_0 S(t-\tau)\ast Q[u_0](s) ds,
\end{equation}
where $S(t,x)$ is the fundamental solution of the linear convolution equation:
\[
\begin{cases}
S_t -D\int_{-\infty}^\infty J(x-y)[S(y,t)-S(x,t)] dy=0,\\
S(0,x)=\delta(x), \mbox{ the Delta function},
\end{cases}
\]
with an explicit form of
\[
S(t,x)=e^{-Dt} \delta(x) +K(t,x),
\]
where $K(t,x)$ is a smooth function defined in Fourier variables by
\[
\hat{K}(t,\xi)=e^{-Dt}(e^{D\hat{J}(\xi)t}-1).
\]
Let us define an iteration to \eqref{new-1}:
\begin{equation}\label{new-2}
u^{(n+1)}(t,x)=S(t)\ast u_0 - \int^t_0 S(t-\tau)\ast P[u^{(n)}] (s) ds + \int^t_0 S(t-\tau)\ast Q[u_0](s) ds.
\end{equation}
When $u^{(n)}\in C_\text{unif}[0,r]$, it is easy to see $u^{(n+1)}\in C_\text{unif}[0,r]$. Since $P[u]$ is Fr\'echet differentiable with respect to $u$
and $DP[u]=d'(\phi+u)$ is a bounded and positive operator, by the existence and uniqueness theory for the convolution equations \cite{Rossi}, we
then prove that such an iteration is a Cauchy sequence with a unique limit:
\[
\lim_{n\to \infty} u^{(n)}(t,x)=u(t,x), \ \mbox{ in } C_{unif}[-r,r],
\]
in another word, the solution for \eqref{new-1} uniquely exists in $C_\text{unif}[-r,r]$.

Next step is  to consider \eqref{eq-pt} for $t\in[r,2r]$.
Since $t-r\in[0,r]$ and $u(t-r,\xi)$ has been solved already,
thus $Q[u](t-r,\xi)$ is known function.
As showed before, we can similarly prove the existence and uniqueness of
the solution to \eqref{eq-pt} in  $C_{unif}[r,2r]$, so then in $C_\text{unif}[-r,2r]$
By repeating this procedure for $t\in[nr,(n+1)r]$ with $n\in\mathbb Z^+$,
we can prove that the perturbed problem \eqref{eq-pt} admits one unique global solution
$u(t,\xi)$ in $C_\text{unif}[-r,+\infty)$.
$\hfill\Box$

\begin{lemma} \label{le-convergence}
There exist a large number $\xi_0\in\mathbb{R}$ and constants $\mu_1>0$, $C>0$,
such that
\begin{equation} \label{eq-convergence}
\|u(t)\|_{L^\infty([\xi_0,+\infty))}\le
Ce^{-\mu_1t}\|u_0\|_{L^\infty([-r,0]\times\mathbb{R})}.
\end{equation}
\end{lemma}
{\it\bfseries Proof.}
This proof is similar to \cite{HuangDCDS16}.
Here we omit it.
$\hfill\Box$

Since $u_-=0$ is the unstable node of \eqref{eq-pt}, heuristically,
for a general initial data $u_0$, we cannot expect the convergence $u\to0$ as $t\to\infty$.
Inspired by \cite{Chern}, we expect the solution $u$ decay to zero
when the initial perturbation is exponentially decay at the far field $\xi=-\infty$.
Thus let us define
\begin{equation} \label{eq-utilde}
\tilde u(t,\xi):=|w(\xi)|^{1/2}u(t,\xi)=e^{-\lambda\xi}u(t,\xi),
\end{equation}
where $\lambda\in(\lambda_1,\lambda_2)$ for $c>c^*$
and $\lambda=\lambda^*$ for $c=c^*$.
Then we substitute $u(t,\xi)=e^{\lambda\xi}\tilde u(t,\xi)$ into \eqref{eq-pt}
and derive the following problem
\begin{equation} \label{eq-ptwt}
\begin{cases}
\displaystyle
\pd{\tilde u}{t}
+c\pd{\tilde u}{\xi}+c\lambda \tilde u
-D\Big(\int_{\mathbb R}J(y)e^{-\lambda y}\tilde u(t,\xi-y)dy-\tilde u(t,\xi)\Big)
+e^{-\lambda \xi}P[e^{\lambda\xi}\tilde u](t,\xi) \\
\displaystyle
\qquad\qquad=e^{-\lambda\xi}\int_{\mathbb R}K(y)Q[e^{\lambda\xi}\tilde u](t-r,\xi-y-cr) dy,
\quad \xi\in\mathbb R,~t>0\\
\tilde u(s,\xi)=e^{-\lambda\xi}u_0(s,\xi)=:\tilde u_0(s,\xi),
\quad \xi\in\mathbb R, ~s\in[-r,0],
\end{cases}
\end{equation}

In order to derive the boundedness of oscillatory traveling waves,
we compare it with the following linear delayed nonlocal dispersion equation
\begin{equation} \label{eq-linear}
\begin{cases}
\displaystyle
\pd{u^+}{t}+c\pd{u^+}{\xi}
-D\int_{\mathbb R}J(y)e^{-\lambda y}u^+(t,\xi-y)dy
+(d'(0)+c\lambda+D)u^+(t,\xi) \\
\displaystyle
\qquad\qquad=b'(0)\int_{\mathbb R}K(y)e^{-\lambda(y+cr)}u^+(t-r,\xi-y-cr) dy,
\quad \xi\in\mathbb R,~t>0\\
u^+(s,\xi)=u_0^+(s,\xi)\ge0,
\quad \xi\in\mathbb R, ~s\in[-r,0],
\end{cases}
\end{equation}

\begin{lemma} \label{le-positive}
For non-negative initial value $u_0^+(s,\xi)\ge0$ on $[-r,0]\times \mathbb R$,
\eqref{eq-linear} has a unique solution satisfying
$u^+(t,\xi)\ge0$ for all $(t,\xi)\in[-r,+\infty)\times \mathbb R$.
\end{lemma}
{\it\bfseries Proof.}
For $t\in[0,r]$, we see that $t-r\in[-r,0]$, and
\begin{align*}
&\pd{u^+}{t}+c\pd{u^+}{\xi}
-D\int_{\mathbb R}J(y)e^{-\lambda y}u^+(t,\xi-y)dy
+(d'(0)+c\lambda+D)u^+(t,\xi) \\
&=b'(0)\int_{\mathbb R}K(y)e^{-\lambda(y+cr)}u^+(t-r,\xi-y-cr) dy\\
&=b'(0)\int_{\mathbb R}K(y)e^{-\lambda(y+cr)}u_0^+(t-r,\xi-y-cr) dy\ge0,
\quad \xi\in\mathbb R,~t>0.
\end{align*}
Making change of variable, $u^+(t,\xi)=u^+(t,x+ct)=:\bar u^+(t,x)$,
\begin{align*}
&\pd{\bar u^+}{t}
-D\int_{\mathbb R}J(y)e^{-\lambda y}\bar u^+(t,x-y)dy
+(d'(0)+c\lambda+D)\bar u^+(t,x)\ge0,
\quad x\in\mathbb R,~t>0.
\end{align*}
As showed in the textbook \cite{Rossi} for nonlocal diffusion equations,  the linear convolution equation \eqref{eq-linear}
exists a unique solution satisfying $\bar u^+(t,x) \ge 0$ for all $x\in\mathbb R$ and $t\in[0,r]$.
We can complete this proof step by step for $t\in[nr,(n+1)r]$ with $n\in\mathbb Z^+$.
$\hfill\Box$

\begin{lemma} \label{le-bounded}
Let $\tilde u(t,\xi)$ and $u^+(t,\xi)$ be the solutions of \eqref{eq-ptwt} and
\eqref{eq-linear}, respectively.
Then
$$|\tilde u(t,\xi)|\le u^+(t,\xi), \qquad (t,\xi)\in[0,+\infty)\times \mathbb R,$$
provided that the initial value
$$|\tilde u_0(t,\xi)|\le u_0^+(t,\xi), \qquad (t,\xi)\in[-r,0]\times \mathbb R.$$
\end{lemma}
{\it\bfseries Proof.}
For $t\in[0,r]$, let
$$U(t,\xi):=u^+(t,\xi)-\tilde u(t,\xi).$$
Since $t-r\in[-r,0]$, according to the initial condition,
we have
$$|\tilde u(t-r,\xi)|=|\tilde u_0(t-r,\xi)|
\le u_0^+(t-r,\xi)=u^+(t-r,\xi),
\quad t\in[0,r], ~\xi\in\mathbb R.$$
Then $U(t,\xi)$ satisfies,
\begin{align*}
&\pd{U}{t}+c\pd{U}{\xi}
-D\int_{\mathbb R}J(y)e^{-\lambda y}U(t,\xi-y)dy
+(d'(0)+c\lambda+D)U(t,\xi) \\
=&b'(0)\int_{\mathbb R}K(y)e^{-\lambda(y+cr)}u^+(t-r,\xi-y-cr) dy\\
&\qquad-e^{-\lambda\xi}\int_{\mathbb R}K(y)Q[e^{\lambda\xi}\tilde u](t-r,\xi-y-cr) dy\\
&\qquad\qquad+e^{-\lambda \xi}P[e^{\lambda\xi}\tilde u](t,\xi)-d'(0)\tilde u(t,\xi)\\
\ge&b'(0)\int_{\mathbb R}K(y)e^{-\lambda(y+cr)}u^+(t-r,\xi-y-cr) dy\\
&\qquad-e^{-\lambda\xi}\int_{\mathbb R}K(y)b'(0)e^{\lambda(\xi-y-cr)}\tilde u(t-r,\xi-y-cr) dy\\
&\qquad\qquad+e^{-\lambda \xi}d'(0)e^{\lambda\xi}\tilde u(t,\xi)-d'(0)\tilde u(t,\xi)\\
=&0, \qquad \xi\in\mathbb R,~t>0,
\end{align*}
and the initial condition
$$U(s,\xi)=u_0^+(s,\xi)-\tilde u(s,\xi)\ge0, \qquad \xi\in\mathbb R, ~s\in[-r,0].$$
Similar to Lemma \ref{le-positive}, we can prove that
$U(t,\xi)\ge0$ (i.e. $u^+(t,\xi)\ge\tilde u(t,\xi)$)
for all $t\in[0,r]$ and $\xi\in\mathbb R$.
By replacing $U(t,\xi):=u^+(t,\xi)-\tilde u(t,\xi)$
with $U(t,\xi):=u^+(t,\xi)+\tilde u(t,\xi)$,
we can similarly prove that $u^+(t,\xi)+\tilde u(t,\xi)\ge0$
for all $t\in[0,r]$ and $\xi\in\mathbb R$.
That is,
\begin{equation} \label{eq-zU}
|\tilde u(t,\xi)|\le u^+(t,\xi), \qquad (t,\xi)\in[0,r]\times \mathbb R.
\end{equation}
Next, when $t\in[r,2r]$, namely, $t-r\in[0,r]$, based on \eqref{eq-zU}
we can similarly prove
\begin{equation*}
|\tilde u(t,\xi)|\le u^+(t,\xi), \qquad (t,\xi)\in[r,2r]\times \mathbb R.
\end{equation*}
Repeating this procedure, we further complete this lemma.
$\hfill\Box$

Now let us recall the fundamental solution theory for time-delayed equations.

\begin{lemma}[\cite{delay}]
Let $z(t)$ be the solution to the following linear time-delayed ODE with time-delay $r>0$
and two constants $k_1$ and $k_2$
\begin{equation} \label{eq-delayODE}
\begin{cases}
\displaystyle
\frac{d}{dt}z(t)+k_1z(t)=k_2z(t-r),\\
z(s)=z_0(s), \quad s\in[-r,0].
\end{cases}
\end{equation}
Then
\begin{equation} \label{eq-delayODEs}
z(t)=e^{-k_1(t+r)}e_r^{\bar k_2t}z_0(-r)+
\int_{-r}^0e^{-k_1(t-s)}e_r^{\bar k_2(t-r-s)}[z_0'(s)+k_1z_0(s)]ds,
\end{equation}
where $\bar k_2:=k_2e^{k_1r}$, and $e_r^{\bar k_2t}$ is the so-called delayed exponential function
in the form
\begin{equation*}
e_r^{\bar k_2t}=
\begin{cases}
0, \quad & -\infty<t<-r,\\
1, \quad & -r\le t<0, \\
1+\frac{\bar k_2t}{1!}, \quad & 0\le t<r,\\
1+\frac{\bar k_2t}{1!}+\frac{\bar k_2^2(t-r)^2}{2!}, \quad & r\le t<2r, \\
\vdots \quad &\vdots\\
1+\frac{\bar k_2t}{1!}+\frac{\bar k_2^2(t-r)^2}{2!}
+\cdots+\frac{\bar k_2^m[t-(m-1)r]^m}{m!}, \quad & (m-1)r\le t<mr, \\
\vdots \quad & \vdots
\end{cases}
\end{equation*}
and $e_r^{\bar k_2t}$ is the fundamental solution to
\begin{equation} \label{eq-delayODEh}
\begin{cases}
\displaystyle
\frac{d}{dt}z(t)=\bar k_2z(t-r),\\
z(s)\equiv1, \quad s\in[-r,0].
\end{cases}
\end{equation}
\end{lemma}

The property of the solution to the delayed linear ODE \eqref{eq-delayODE}
is well-known \cite{MeiWang2011}.

\begin{lemma}[\cite{MeiWang2011}] \label{le-delayODE}
Let $k_1\ge0$ and $k_2\ge0$.
Then the solution $z(t)$ to \eqref{eq-delayODE}
(or equivalently the function $z(t)$ in \eqref{eq-delayODEs}) satisfies
$$
|z(t)|\le C_0e^{-k_1t}e_r^{\bar k_2t},
$$
where
$$
C_0:=e^{-k_1r}|z_0(-r)|+\int_{-r}^0e^{k_1s}|z_0'(s)+k_1z_0(s)|ds,
$$
and the fundamental solution $e_r^{\bar k_2t}$ with $\bar k_2>0$ to \eqref{eq-delayODEh}
satisfies
$$
e_r^{\bar k_2t}\le C(1+t)^{-\gamma}e^{\bar k_2t}, \quad t>0,
$$
for arbitrary number $\gamma>0$.
Furthermore, when $k_1\ge k_2\ge0$, there exists a constant $0<\varepsilon_1<1$ such that
$$
e^{-k_1t}e_r^{\bar k_2t}\le Ce^{-\varepsilon(k_1-k_2)t}, \quad t>0,
$$
and the solution $z(t)$ to \eqref{eq-delayODE} satisfies
$$
|z(t)|\le Ce^{-\varepsilon(k_1-k_2)t}, \quad t>0.
$$
\end{lemma}

For the linear delayed nonlocal dispersion equation \eqref{eq-linear},
we take Fourier transform of $u^+(t,\xi)$ and denote it by $\hat u^+(t,\eta)$
or $\mathscr{F}[u^+](t,\eta)$,
then we have
\begin{equation} \label{eq-fourier}
\begin{cases}\displaystyle
\frac{d}{dt}\hat u^+(t,\eta)+A(\eta)\hat u^+(t,\eta)=B(\eta)\hat u^+(t-r,\eta), \\
\hat u^+(s,\eta)=\hat u_0^+(s,\eta), \quad s\in[-r,0], ~\eta\in\mathbb R,
\end{cases}
\end{equation}
where
\begin{equation} \label{eq-AB}
\begin{cases}
A(\eta):=-D\mathscr{F}[J(y)e^{-\lambda y}](\eta)+(d'(0)+c\lambda+D)+ic\eta,\\
B(\eta):=b'(0)e^{icr\eta}\mathscr{F}[K(y)e^{-\lambda(y+cr)}](\eta),
\end{cases}
\quad c\ge c^*.
\end{equation}

The liner delayed equation \eqref{eq-fourier} can be solved by
\begin{align} \nonumber
\hat u^+(t,\eta)=&e^{-A(\eta)(t+r)}e_r^{\bar B(\eta)t}\hat u_0^+(-r,\eta)
\\ \label{eq-solution}
&+\int_{-r}^0e^{-A(\eta)(t-s)}e_r^{\bar B(\eta)(t-r-s)}
[\frac{d}{ds}u_0^+(s,\eta)+A(\eta)u_0^+(s,\eta)]ds,
\end{align}
where $\bar B(\eta):=B(\eta)e^{A(\eta)r}$.
Then by taking the inverse Fourier transform, we get
\begin{align} \nonumber
u^+(t,x)=&\frac{1}{2\pi}\int_\mathbb{R}
e^{ix\cdot\eta}e^{-A(\eta)(t+r)}e_r^{\bar B(\eta)t}\hat u_0^+(-r,\eta)d\eta
\\ \label{eq-solutionF}
&+\int_{-r}^0\frac{1}{2\pi}\int_\mathbb{R}e^{ix\cdot\eta}
e^{-A(\eta)(t-s)}e_r^{\bar B(\eta)(t-r-s)}
[\frac{d}{ds}u_0^+(s,\eta)+A(\eta)u_0^+(s,\eta)]d\eta ds,
\end{align}
and its derivatives
\begin{align} \nonumber
\partial_x^ku^+(t,x)=&\frac{1}{2\pi}\int_\mathbb{R}
e^{ix\cdot\eta}(i\eta)^ke^{-A(\eta)(t+r)}e_r^{\bar B(\eta)t}\hat u_0^+(-r,\eta)d\eta
\\ \label{eq-solutionD}
&+\int_{-r}^0\frac{1}{2\pi}\int_\mathbb{R}e^{ix\cdot\eta}(i\eta)^k
e^{-A(\eta)(t-s)}e_r^{\bar B(\eta)(t-r-s)}
[\frac{d}{ds}u_0^+(s,\eta)+A(\eta)u_0^+(s,\eta)]d\eta ds,
\end{align}
for $k\in\mathbb Z^+$.
For simplicity, we denote
\begin{align} \label{eq-I1}
I_1(t,\eta)&:=(i\eta)^ke^{-A(\eta)(t+r)}e_r^{\bar B(\eta)t}\hat u_0^+(-r,\eta),
\\ \label{eq-I2}
I_2(t-s,\eta)&:=(i\eta)^ke^{-A(\eta)(t-s)}e_r^{\bar B(\eta)(t-r-s)}
[\frac{d}{ds}u_0^+(s,\eta)+A(\eta)u_0^+(s,\eta)].
\end{align}

Next, we are going to estimate the decay rates for the solution $u^+(t,x)$.

\begin{lemma} \label{le-Parseval}
Suppose that $u_0^+\in C([-r,0];H^{m}(\mathbb R)\cap L^1(\mathbb R))$
and $\partial_s u_0^+\in C([-r,0]; H^{m}(\mathbb R) \cap L^1(\mathbb R))$ for $m\ge1$.
Then there exist constants $C>0$ and $\varepsilon_1>0$ such that
$$\|\partial_x^ku^+(t,x)\|_{L^2(\mathbb R)}
\le C\mathcal{E}^ke^{-\varepsilon_1\mu_0(c)t}t^{-\frac{1+2k}{2\alpha}}, \quad t>0,$$
for $k=0,1,\cdots,[m]$,
where $\mu_0(c)=G_c(\lambda)-H_c(\lambda)>0$ for $c>c^*$
and $\mu_0(c)=0$ for $c=c^*$, and
$$\mathcal{E}^k=\|u_0^+(-r)\|_{L^1(\mathbb R)}+\|u_0^+(-r)\|_{H^k(\mathbb R)}
+\int_{-r}^0\Big[\|(u_0^+,\pd{}{s}u_0^+)(s)\|_{L^1(\mathbb R)}
+\|(u_0^+,\pd{}{s}u_0^+)(s)\|_{H^k(\mathbb R)}\Big]ds.$$
Furthermore,
$$\|u^+(t,x)\|_{L^\infty(\mathbb R)}\le C
\mathcal{E}^1e^{-\varepsilon_1\mu_0(c)t}t^{-\frac{1}{\alpha}}, \quad t>0.$$
\end{lemma}
{\it\bfseries Proof.}
By using Parseval's equality, we have
\begin{align} \nonumber
&\|\partial_x^ku^+(t,x)\|_{L^2(\mathbb R)}=\Big\|\mathscr{F}^{-1}[I_1](t,x)+
\int_{-r}^0\mathscr{F}^{-1}[I_2](t-s,x)ds\Big\|_{L^2(\mathbb R)}\\ \nonumber
\le &\|\mathscr{F}^{-1}[I_1](t,x)\|_{L^2(\mathbb R)}+
\int_{-r}^0\|\mathscr{F}^{-1}[I_2](t-s,x)\|_{L^2(\mathbb R)}ds
\\ \label{eq-Parseval}
=&\|I_1(t,\eta)\|_{L^2(\mathbb R)}+
\int_{-r}^0\|I_2(t-s,\eta)\|_{L^2(\mathbb R)}ds
\end{align}
We note that
\begin{align*}
|e^{-A(\eta)t}|
&=e^{-(d'(0)+c\lambda+D)t}\Big|\exp\Big(tD\int_{\mathbb R}
J(y)e^{-\lambda y}e^{-iy\cdot\eta}dy\Big)\Big|\\
&=e^{-(d'(0)+c\lambda+D)t}\exp\Big(tD\int_{\mathbb R}
J(y)e^{-\lambda y}\cos(y\cdot\eta)dy\Big)\\
&=\exp\Big(-t\big(d'(0)+c\lambda+D-D\int_\mathbb{R}J(y)e^{-\lambda y}dy\big)\Big)\\
&\qquad \cdot\exp\Big(-tD\int_{\mathbb R}
J(y)e^{-\lambda y}(1-\cos(y\cdot\eta))dy\Big),
\end{align*}
and since $J(y)$ is even and $\sin(y\cdot\eta)$ is odd, we have
\begin{align*}
&\exp\Big(-tD\int_{\mathbb R}J(y)e^{-\lambda y}(1-\cos(y\cdot\eta))dy\Big)\\
=&\exp\Big(-tD\int_{\mathbb R}J(y)\frac{e^{-\lambda y}+e^{\lambda y}}{2}
(1-\cos(y\cdot\eta))dy\Big)\\
\le &\exp\Big(-tD\int_{\mathbb R}J(y)(1-\cos(y\cdot\eta))dy\Big)\\
=&\exp\Big(-tD\int_{\mathbb R}J(y)(1-\cos(y\cdot\eta)-i\sin(y\cdot\eta))dy\Big)\\
=&\exp(tD(\hat J(\eta)-1).
\end{align*}
Therefore,
\begin{align*}
|e^{-A(\eta)t}|\le \exp\Big(-t\big(d'(0)+c\lambda+D-D\int_\mathbb{R}J(y)e^{-\lambda y}dy\big)\Big)\cdot\exp(tD(\hat J(\eta)-1)=:e^{-k_1t},
\end{align*}
with
$$k_1=d'(0)+c\lambda+D-D\int_\mathbb{R}J(y)e^{-\lambda y}dy
+D(1-\hat J(\eta)).$$
Also, we have
\begin{align*}
|B(\eta)|=b'(0)|\mathscr{F}[K(y)e^{-\lambda(y+cr)}](\eta)|
\le b'(0)\int_\mathbb{R}K(y)e^{-\lambda(y+cr)}dy=:k_2,
\end{align*}
According to the assumption of the existence of traveling waves, there holds
\begin{equation} \label{eq-GH}
G_c(\lambda):=c\lambda-D\int_{\mathbb R}J(y)e^{-\lambda y}dy+D+d'(0)
\ge b'(0)\int_{\mathbb R}K(y)e^{-\lambda(y+cr)}dy=:H_c(\lambda),
\end{equation}
for $\lambda\in(\lambda_1,\lambda_2)$ with $c>c^*$
or $\lambda=\lambda^*$ with $c=c^*$.
Therefore,
\begin{align*}
k_1-k_2\ge G_c(\lambda)-H_c(\lambda)+D(1-\hat J(\eta))\ge
\mu_0(c)+D(1-\hat J(\eta)),
\end{align*}
where $\mu_0(c)=G_c(\lambda)-H_c(\lambda)>0$ for $c>c^*$
and $\mu_0(c)=0$ for $c=c^*$.

From the assumption, there exist constants $0<\kappa_1\le\kappa_2$, $0<\delta<1$,
and $\tilde r>0$, such that
\begin{equation} \label{eq-J}
\begin{cases}
\kappa_1|\eta|^\alpha\le 1-\hat J(\eta)\le \kappa_2|\eta|^\alpha,
\quad &|\eta|\le\tilde r,\\
\delta\le 1-\hat J(\eta)\le 1, \quad &|\eta|>\tilde r.
\end{cases}
\end{equation}

Using the above estimates in Lemma \ref{le-delayODE} for time-delayed ODE, we obtain
\begin{align*}
\|I_1(t,\eta)\|_{L^2(\mathbb R)}^2
=&\int_\mathbb{R}|e^{-A(\eta)(t+r)}e_r^{\bar B(\eta)t}|^2
|\eta|^{2k}|\hat u_0^+(-r,\eta)|^2d\eta\\
\le &C\int_\mathbb{R}(e^{-k_1(t+r)}e_r^{\bar k_2t})^2
|\eta|^{2k}|\hat u_0^+(-r,\eta)|^2d\eta\\
\le &C\int_\mathbb{R}(e^{-\varepsilon_1(k_1-k_2)t})^2
|\eta|^{2k}|\hat u_0^+(-r,\eta)|^2d\eta\\
\le &Ce^{-2\varepsilon_1\mu_0(c)t}\int_\mathbb{R}e^{-2\varepsilon_1D(1-\hat J(\eta))t}
|\eta|^{2k}|\hat u_0^+(-r,\eta)|^2d\eta,
\end{align*}
and furthermore
\begin{align*}
&\int_\mathbb{R}e^{-2\varepsilon_1D(1-\hat J(\eta))t}
|\eta|^{2k}|\hat u_0^+(-r,\eta)|^2d\eta \\
\le&\int_{|\eta|\le\tilde r}e^{-2\varepsilon_1D(1-\hat J(\eta))t}
|\eta|^{2k}|\hat u_0^+(-r,\eta)|^2d\eta
+\int_{|\eta|>\tilde r}e^{-2\varepsilon_1D(1-\hat J(\eta))t}
|\eta|^{2k}|\hat u_0^+(-r,\eta)|^2d\eta\\
\le&\int_{|\eta|\le\tilde r}e^{-2\varepsilon_1D\kappa_1|\eta|^\alpha t}
|\eta|^{2k}|\hat u_0^+(-r,\eta)|^2d\eta
+\int_{|\eta|>\tilde r}e^{-2\varepsilon_1D\delta t}
|\eta|^{2k}|\hat u_0^+(-r,\eta)|^2d\eta
\\
\le &\|\hat u_0^+(-r,\eta)\|_{L^\infty(\mathbb R)}^2t^{-\frac{1+2k}{\alpha}}
\int_{|\eta|\le\tilde r}e^{-2\varepsilon_1D\kappa_1|\eta t^{1/\alpha}|^\alpha}
|\eta t^{1/\alpha}|^{2k}d(\eta t^{1/\alpha})\\
&\qquad+e^{-2\varepsilon_1D\delta t}\int_{|\eta|>\tilde r}
|\eta|^{2k}|\hat u_0^+(-r,\eta)|^2d\eta\\
\le &C(\|u_0^+(-r,x)\|_{L^1(\mathbb R)}^2+\|u_0^+(-r,x)\|_{H^k(\mathbb R)}^2)
t^{-\frac{1+2k}{\alpha}}.
\end{align*}
Substitute it into the above inequality, we obtain
\begin{align*}
\|I_1(t,\eta)\|_{L^2(\mathbb R)}
\le C(\|u_0^+(-r,x)\|_{L^1(\mathbb R)}+\|u_0^+(-r,x)\|_{H^k(\mathbb R)})
e^{-\varepsilon_1\mu_0(c)t}t^{-\frac{1+2k}{2\alpha}}.
\end{align*}

Thus, in a similar way, we can also prove
\begin{align*}
&\|I_2(t-s,\eta)\|_{L^2(\mathbb R)}\\
=&\Big(\int_\mathbb{R}|\eta|^{2k}|e^{-A(\eta)(t-s)}e_r^{\bar B(\eta)(t-r-s)}|^2
\big|\left[\frac{d}{ds}u_0^+(s,\eta)+A(\eta)u_0^+(s,\eta)\right]\big|^2d\eta\Big)^\frac{1}{2}\\
\le& C e^{-\varepsilon_1\mu_0(c)t}\Big(
\int_\mathbb{R}e^{-2\varepsilon_1D(1-\hat J(\eta))t}
|\eta|^{2k}\big|\left[\frac{d}{ds}u_0^+(s,\eta)+A(\eta)u_0^+(s,\eta)\right]\big|^2d\eta
\Big)^\frac{1}{2}\\
\le& C(\|u_0^+(s,x)\|_{L^1(\mathbb R)}
+\|\partial_su_0^+(s,x)\|_{L^1(\mathbb R)}
+\|u_0^+(s,x)\|_{H^k(\mathbb R)}+\|\partial_su_0^+(s,x)\|_{H^k(\mathbb R)})\\
&\quad \cdot e^{-\varepsilon_1\mu_0(c)t}t^{-\frac{1+2k}{2\alpha}}.
\end{align*}
The proof is completed.
$\hfill\Box$

\begin{lemma} \label{le-decay}
When $\tilde u_0\in C([-r,0];H^1(\mathbb R)\cap L^1(\mathbb R))$
and $\partial_s\tilde u_0\in C([-r,0];H^1(\mathbb R)\cap L^1(\mathbb R))$.
Then there exist constants $C>0$ and $\varepsilon_1>0$ such that
$$\|\tilde u(t,x)\|_{L^\infty(\mathbb R)}\le C
e^{-\varepsilon_1\mu_0(c)t}t^{-\frac{1}{\alpha}}, \quad t>0,$$
where $\mu_0(c)=G_c(\lambda)-H_c(\lambda)>0$ for $c>c^*$
and $\mu_0(c)=0$ for $c=c^*$
\end{lemma}
{\it\bfseries Proof.}
We can choose $u_0^+\in C([-r,0];H^1(\mathbb R)\cap L^1(\mathbb R))$
and $\partial_su_0^+\in C([-r,0];H^1(\mathbb R)\cap L^1(\mathbb R))$
such that $|\tilde u_0(s,\xi)\le u_0^+(s,\xi)|$ for all $\xi\in\mathbb R$ and $s\in[-r,0]$.
Combining the boundedness Lemma \ref{le-bounded} and the decay estimate Lemma \ref{le-decay},
we immediately get the convergence result.
$\hfill\Box$

{\it\bfseries Proof of Theorem \ref{th-stability}.}
Based on Lemma \ref{le-convergence} and \ref{le-decay}, we have the
convergence rates of $u(t,x)$.
$\hfill\Box$

\section{Numerical Simulation of traveling waves}

In this section, we numerically study the stability of travelling waves of \eqref{eq-main}
by using the finite-difference method and iteration technique.
We consider the Nicholson's
blowflies equation with nonlocal diffusion by taking the kernels as the heat kernel
\begin{equation}\nonumber
J(x)=\frac{1}{\sqrt{4\pi\alpha}}e^{-x^2/(4\alpha)}=:f_\alpha(x), \quad
K(x)=\frac{1}{\sqrt{4\pi\beta}}e^{-x^2/(4\beta)}=:f_\beta(x).
\end{equation}
And we choose the Nicholson's type of birth rate and death rate
$$
b(v)=pve^{-av},\qquad d(v)=\delta v.
$$
They satisfies the hypotheses (H1)--(H2).
This equation possesses two constant equilibria: $v_-=0$
and $v_+=\ln(p/\delta)/a$.
The traveling wave equation \eqref{eq-tw} now reads
\begin{equation} \label{eq-twNum}
\begin{cases}
c\phi'-D(f_\alpha*\phi-\phi)+\delta\phi=pf_\beta*(\phi_{cr}e^{-a\phi_{cr}}), \\
\phi(-\infty)=0, \quad \phi(+\infty)=v_+,
\end{cases}
\end{equation}
where $\phi_{cr}(\xi)=\phi(\xi-cr)$.
The initial value problem \eqref{eq-main} in the moving coordinates $(t,\xi)$ is
\begin{equation} \label{eq-mainNum}
\begin{cases}
\displaystyle
\pd{v}{t}+c\pd{v}{\xi}-D(f_\alpha*v-v)+\delta v=pf_\beta*(v_re^{-av_r}),\\
v(s,\xi)=v_0(s,\xi), \quad s\in[-r,0], ~ \xi\in\mathbb R,
\end{cases}
\end{equation}
where $v_r(t,\xi)=v(t-r,\xi-cr)$.

The framework for the numerical simulation for local dispersion case can refer to \cite{XJMYJDE}.
Here we present a framework for nonlocal dispersion problem.
The initial value problem \eqref{eq-mainNum} is
required to be solved on $\mathbb R^+\times\mathbb R$,
but numerically we have to impose a finite computational domain $(-M,M)$ for spatial variable $\xi$
and a finite time interval $(0,T)$ with some selected large numbers $M$ and $T$.
Consider the following second order differential problem with
artificial viscosities $-\mu\Delta$ on both sides of the equation
\begin{equation} \label{eq-Num}
\begin{cases}
\displaystyle
\pd{v}{t}+c\pd{v}{\xi}-\mu\frac{\partial^2v}{\partial x^2}
+(D+\delta)v=Df_\alpha*\widetilde v
+pf_\beta*(\widetilde v_re^{-a\widetilde v_r})
-\mu\frac{\partial^2\widetilde v}{\partial x^2}, \\
\qquad (t,x)\in(0,T)\times(-M,M),\\
v(t,-M)=0, \qquad v(t,M)=v_+, \quad t\in(0,T), \\
v(s,\xi)=v_0(s,\xi), \quad s\in[-r,0], ~ \xi\in(-M,M),
\end{cases}
\end{equation}
where $\widetilde v\in C^{1,2}([-r,T]\times[-M,M])$,
$\widetilde v(s,\xi)=v_0(s,\xi)$ for $s\in[-r,0]$, $\xi\in[-M,M]$, and
$\widetilde v(t,-M)=0$, $\widetilde v(t,M)=K$,
$\mu>0$ is the regularization parameter.
The solution of \eqref{eq-Num} is denoted by $v=G(\widetilde v)$.
If the operator $G$ admits a fixed point $\phi$ such that $\phi=G(\phi)$,
we may regard it as an approximate solution of \eqref{eq-mainNum}.
The nonhomogeneous linear problem \eqref{eq-Num} can be solved by the standard
finite-difference method such as the backward difference scheme
for any given $\widetilde v\in C^2([-r,T]\times[-M,M])$.
The reason here we add an artificial viscosities $-\mu\varDelta$ on both sides
is that our simulation without regularization seems not to be stable
as its numerical solution blows up in finite time.

The traveling wave equation \eqref{eq-twNum} can also be numerically simulated as follows.
Consider the following second order ordinary differential problem with
artificial viscosities $-\mu\varDelta$ on both sides
\begin{equation} \label{eq-Numtw}
\begin{cases}
c\phi'-\mu\phi''
+(D+\delta)\phi=Df_\alpha*\widetilde\phi
+pf_\beta*(\widetilde\phi_{cr}e^{-a\widetilde\phi_{cr}})
-\mu\widetilde\phi'', \quad \xi\in(-M,M),\\
\phi(t,-M)=0, \qquad \phi(t,M)=v_+,
\end{cases}
\end{equation}
where $\widetilde\phi\in C^2([-M,M])$,
$\widetilde\phi(t,-M)=0$, $\widetilde\phi(t,M)=K$.
For any given $\widetilde\phi$, \eqref{eq-Numtw} admits an unique solution denoted by
$\phi=H(\widetilde\phi)$.
The fixed point of operator $H$ (if exists) can be regarded as the traveling wave
solution of \eqref{eq-twNum}.
We numerically solve \eqref{eq-Numtw} by starting from the initial profile
$\widetilde\phi(\xi)=v_+e^{\lambda_*\xi}/(1+e^{\lambda_*\xi})$.

For simplicity,
we choose $D=\delta=v_+=\mu=1$ (which means $a=\log(p/\delta)$)
and leave $p,r$ and the initial data $v_0$ free.
We also take $\alpha=\beta=1$ in the kernels $J$ and $K$.
The critical traveling wave speed $c^*$ is uniquely determined by \eqref{eq-characteristic1}-\eqref{eq-characteristic2}. From our stability theorem \ref{th-stability},
let us choose the initial data of \eqref{eq-main} satisfying
\begin{equation*}
\lim_{x\to-\infty}v_0(s,x)=0, \quad
\lim_{x\to\infty}v_0(s,x)=v_+ \quad
\text{uniformly in}\quad s\in[-r,0],
\end{equation*}
and
\begin{equation*}
e^{-\lambda_*x}|v_0(s,x)-\phi(x+c^*s)|\to 0\quad \text{as} \quad x\to-\infty,\quad \text{uniformly in}\quad s\in[-r,0].
\end{equation*}
So we choose the initial data in the moving coordinates $(t,\xi)$ as
\begin{equation*}
v_0(s,\xi)=\phi(\xi)+\varepsilon f_\gamma(\xi), \quad s\in[-r,0], ~\xi\in\mathbb R,
\end{equation*}
where $\phi$ is the traveling wave solution to \eqref{eq-twNum}.
In simulation, we choose $\phi$ being the numerical solution to \eqref{eq-Numtw}
corresponding to critical speed $c^*$, and $\varepsilon=1$, $\gamma=1$.

The results in \cite{HuangDCDS16,Tang} shows that
when $d'(v_+)<|b'(v_+)|$,
the traveling wave exist for $0<r<\overline r$, and no traveling waves exist for
$r\ge\overline r$, where
\begin{equation}\label{eq-overline-r}
\overline r:=\frac{\pi-\arctan(\sqrt{|b'(v_+)|^2-d'(v_+)^2}/d'(v_+))}
{\sqrt{|b'(v_+)|^2-d'(v+)^2}},
\end{equation}
and when $d'(v_+)\ge|b'(v_+)|$,  the traveling wave globally exists in time for any time delay $r$.
Moreover, the traveling waves are monotone for $0<r<\underline r$
and it may be oscillatory for $r\ge\underline r$,
where $\underline r$ is given by
\begin{equation}\label{eq-underline-r}
|b'(v_+)|\underline r e^{d\underline r+1}=1.
\end{equation}
The condition $d'(v_+)\ge|b'(v_+)|$ is equivalent to $e < p/\delta\le e^2$,
and $d'(v_+)<|b'(v_+)|$ is equivalent to $p/\delta>e^2$.

Next, we report the results in four cases, see Table \ref{tb-pr} for the details.
{\small
\begin{table}[!ht] \label{tb-pr}
\caption{Different cases for selection of $p$, $r$}
\begin{tabular}{lrrlllll}
\hline
Case & $p$ & $r$ & Zone of $p/\delta$ & Zone of $r$
&$c^*$&$\lambda_*$&Behavior of $v$\\
\hline
1   &5 &0.2&$\frac{p}{d}\in(e,e^2)$& $r<\underline r (0.403...)$
&5.104...&0.780...&monotone \\
  2 &5 &2&$\frac{p}{d}\in(e,e^2)$& $r>\underline r$
  &1.310...&0.754...&oscillatory\\
 3  &10 &0.2&$\frac{p}{d}>e^2$& $r<\underline r (0.225...)$
 &7.153...&0.931...&monotone\\
  4 &10  &2&$\frac{p}{d}>e^2$& $\underline r<r<\overline r(2.930...)$
  &1.617...&0.858...&oscillatory\\
   \hline
\end{tabular}
\end{table}
}

\begin{figure}[htb]
\includegraphics[width=0.9\textwidth,height=0.3\textheight]{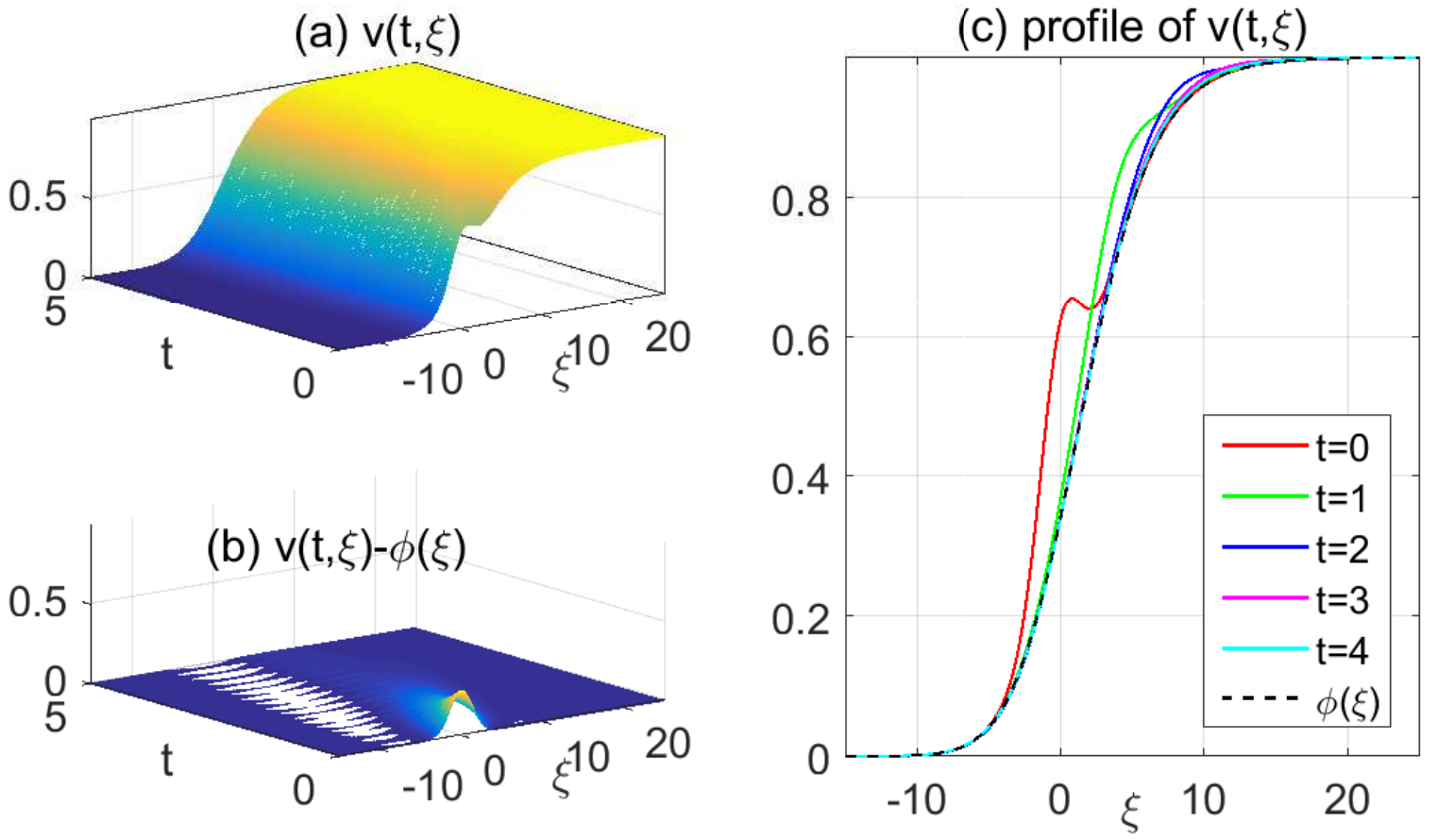}
\caption{Case 1. $e<\frac{p}{d}\le e^2$ with small time delay $r<\underline r$.
(a) 3D-graphs of $v(t, \xi)$;
(b) 3D-graphs of the error $v(t, \xi)-\phi(\xi)$;
and (c) 2D-graphs of $v(t, \xi)$ at
$t = 0, 1, 2, 3, 4$ and the traveling wave $\phi(\xi)$.}
\label{fig-1}
\end{figure}

{\bf Case 1. $e<\frac{p}{d}\le e^2$ and $r<\underline r$ the solution converges to a monotone critical travelling wave $\phi(x+c^*t)$.}
We take $p = 5$ and $r = 0.2$. In this case,  when $e<\frac{p}{d}\le e^2$,
the birth rate function $b(v)$ for $v \in [0, v_+]$ is non-monotone,
where $v_+=\frac{1}{a}\ln\frac{p}{d}$.
A direct calculation from  \eqref{eq-characteristic1}, \eqref{eq-characteristic2},
and \eqref{eq-underline-r} gives
$\underline r=0.4032979$, $c^*=5.1041202$, $\lambda_*=0.7801950$.
Since $r<\underline r$,  the critical wave $\phi(x+c^*t)$ is monotone.
As numerically demonstrated in Figure \ref{fig-1}, we can see that the solution behaves
exactly like a monotone traveling wave,
which are consistent with our stability Theorem \ref{th-stability}.

\begin{figure}[htb]
\centering
\includegraphics[width=0.9\textwidth,height=0.3\textheight]{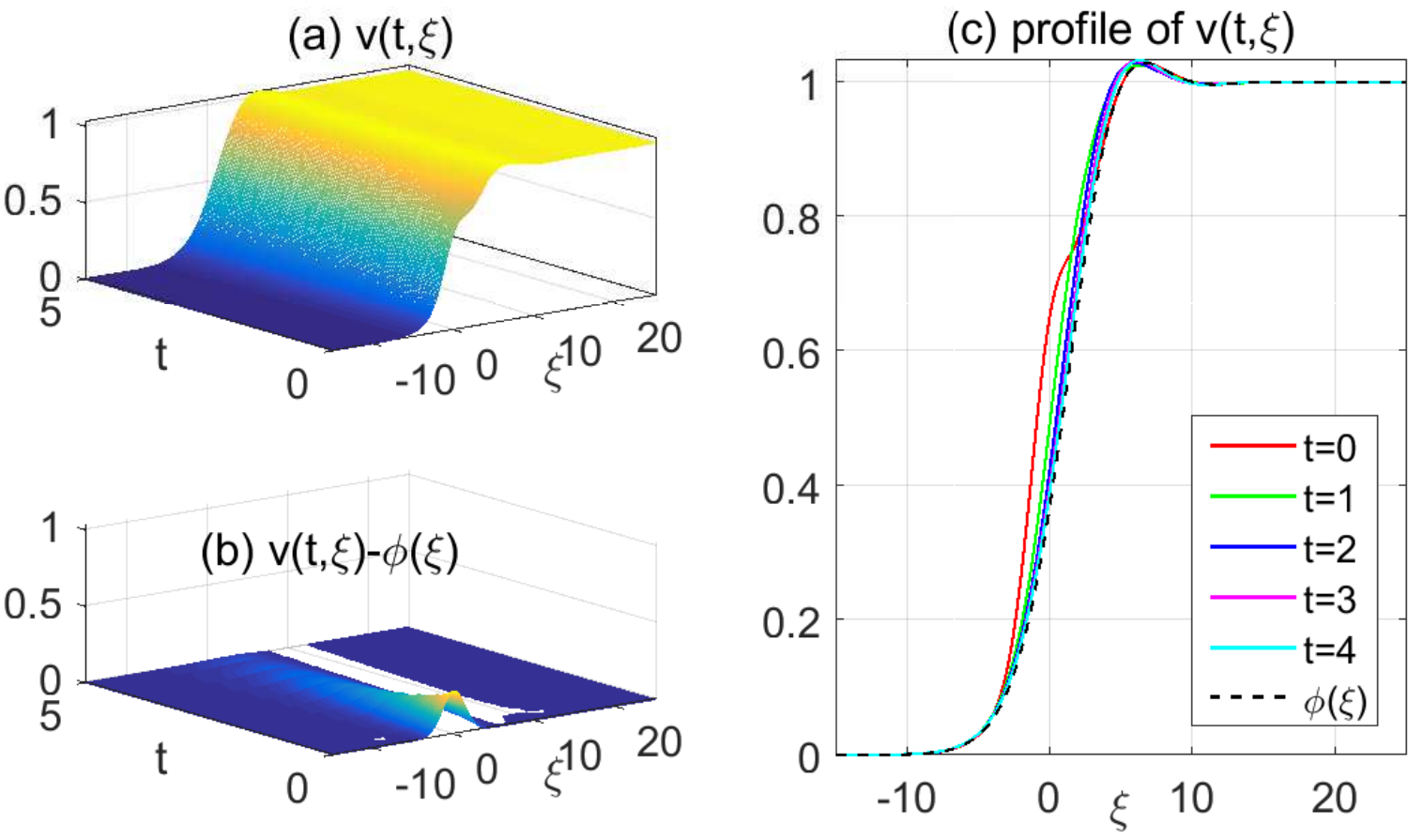}
\caption{Case 2. $e<\frac{p}{d}\le e^2$ with big time delay $r>\underline r$.
(a) 3D-graphs of $v(t, \xi)$;
(b) 3D-graphs of the error $v(t, \xi)-\phi(\xi)$;
and (c) 2D-graphs of $v(t, \xi)$ at
$t = 0, 1, 2, 3, 4$ and the traveling wave $\phi(\xi)$.}
\label{fig-2}
\end{figure}
{\bf Case 2. $e<\frac{p}{d}\le e^2$ and $r>\underline r$ the solution converges to an oscillatory critical travelling wave $\phi(x+c^*t)$.}
We can choose $p=5$ and $r=2$.
Similarly,
we have $\underline r=0.4032979$, $c^*=1.3108958$, $\lambda_*=0.7548876$
from \eqref{eq-characteristic1}, \eqref{eq-characteristic2}, and \eqref{eq-underline-r}.
In this case, $r>\underline r$, the critical traveling wave may be oscillating.
The numerical results are shown in Figure \ref{fig-2}.
Figure \ref{fig-2} shows that the solution $v(t, x)$ is oscillating,
and it converges to the oscillatory critical traveling wave.

\begin{figure}[htb]
\centering
\includegraphics[width=0.9\textwidth,height=0.3\textheight]{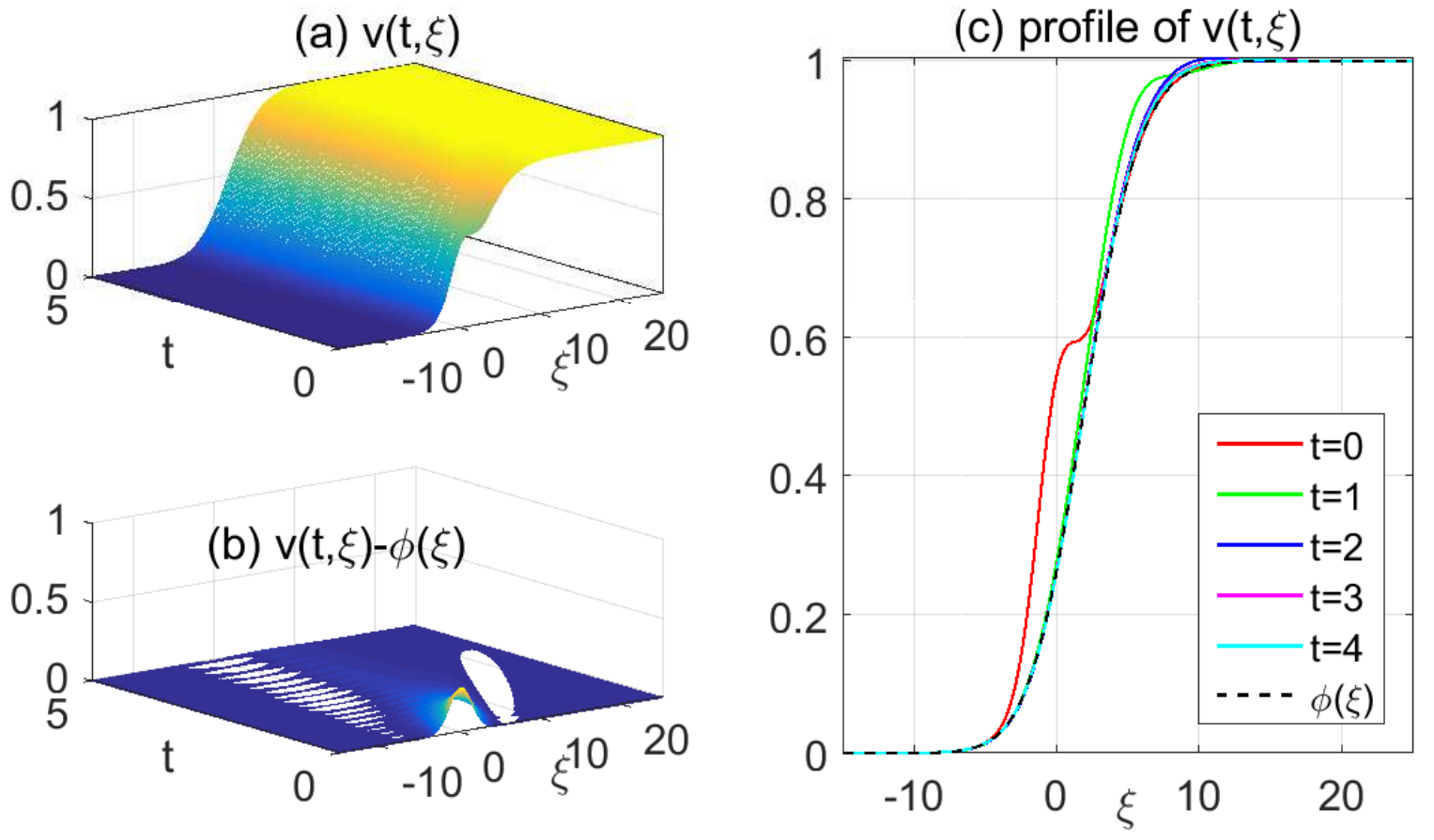}
\caption{Case 3. $\frac{p}{d}> e^2$ with small time delay $r<\underline r$.
(a) 3D-graphs of $v(t, \xi)$;
(b) 3D-graphs of the error $v(t, \xi)-\phi(\xi)$;
and (c) 2D-graphs of $v(t, \xi)$ at
$t = 0, 1, 2, 3, 4$ and the traveling wave $\phi(\xi)$.}
\label{fig-3}
\end{figure}
{\bf Case 3. $\frac{p}{d}> e^2$ and $r<\underline r$ the solution converges to a monotone critical travelling wave $\phi(x+c^*t)$.}
We take $p = 10$ and $r = 0.2$ and get
$\overline r=2.9304424$, $\underline r=0.2254235$, $c^*=7.1531405$, $\lambda_*=0.9315197$
from \eqref{eq-characteristic1}, \eqref{eq-characteristic2},
\eqref{eq-overline-r} and \eqref{eq-underline-r}.
Since $r<\underline r$,  the critical wave $\phi(x+c^*t)$ is monotone.
The numerical results showed in Figure \ref{fig-3} demonstrates that $v(t,x)$
behaves like the monotone critical traveling wave.

\begin{figure}[htb]
\centering
\includegraphics[width=0.9\textwidth,height=0.3\textheight]{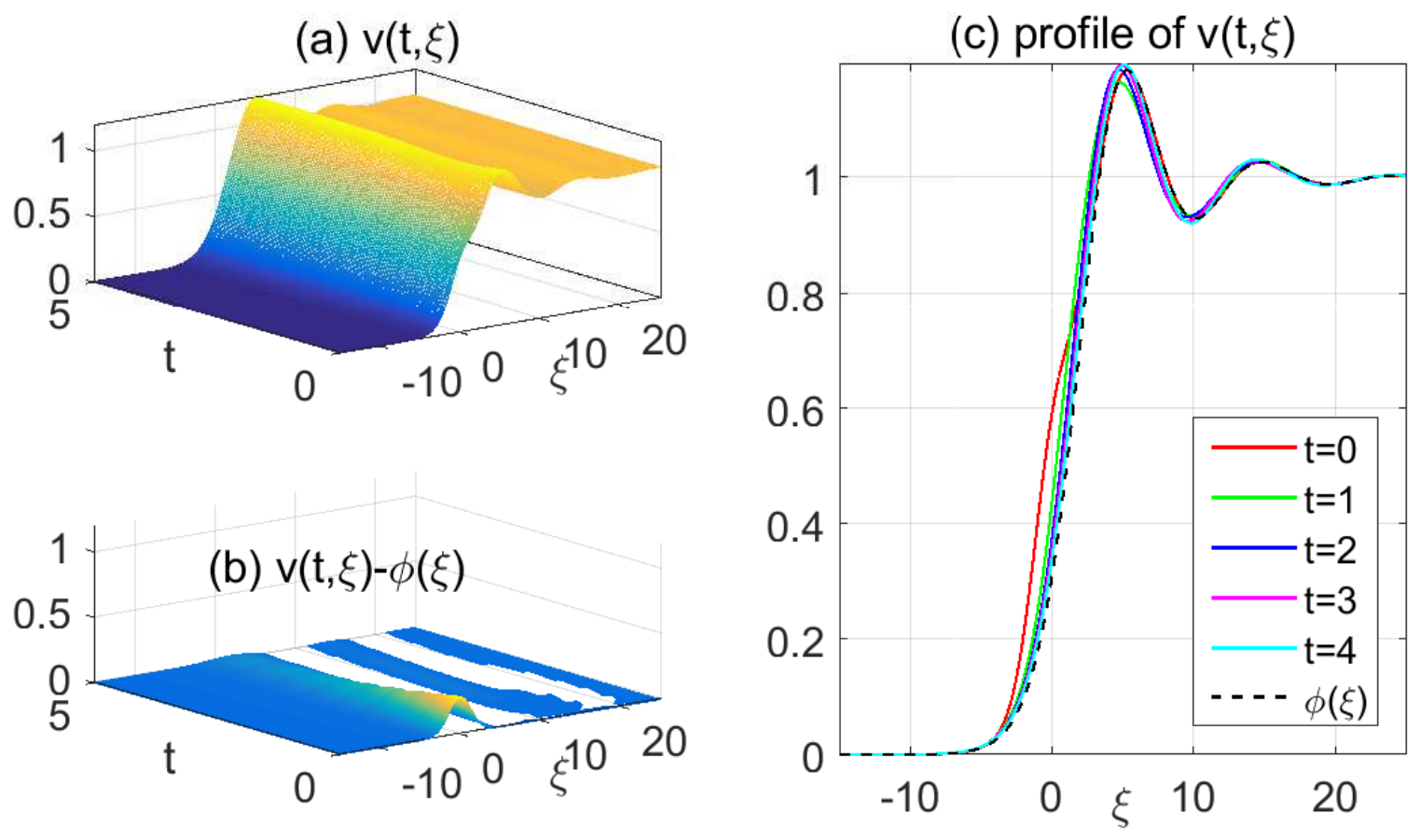}
\caption{Case 4. $\frac{p}{d}> e^2$ with  time delay $\underline r<r<\overline r$.
(a) 3D-graphs of $v(t, \xi)$;
(b) 3D-graphs of the error $v(t, \xi)-\phi(\xi)$;
and (c) 2D-graphs of $v(t, \xi)$ at
$t = 0, 1, 2, 3, 4$ and the traveling wave $\phi(\xi)$.}
\label{fig-4}
\end{figure}
{\bf Case 4. $\frac{p}{d}> e^2$ and $\underline r<r<\overline r$ the solution converges to an oscillatory critical travelling wave $\phi(x+c^*t)$.}
We take $p = 10$ and $r = 2$.
A simple calculation from \eqref{eq-characteristic1}, \eqref{eq-characteristic2},
\eqref{eq-overline-r} and \eqref{eq-underline-r} gives
$\overline r=2.9304424$, $\underline r=0.2254235$, $c^*=1.6178475$, $\lambda_*=0.8586847$.
The numerical results given in Figure \ref{fig-4}.

{\bf Acknowledgement}.
The research of S. Ji was
supported by NSFC Grant No.~11701184,
the Fundamental Research Funds for the Central Universities (No.~2017BQ109),
and the China Postdoctoral Science Foundation (No.~2017M610517).
The research of R. Huang was supported in part
by NSFC Grants No. 11671155 and No. 11771155, NSF of Guangdong Grant No. 2016A030313418,
and NSF of Guangzhou Grant No. 201607010207.
The research of M. Mei was supported in part
by NSERC Grant RGPIN 354724-16, and FRQNT Grant No. 2019-CO-256440.
The research of J. Yin was supported in
part by NSFC Grant No. 11771156.

\end{document}